\documentclass[12pt]{amsart}
\usepackage[colorlinks=true, pdfstartview=FitV, linkcolor=blue, citecolor=blue, urlcolor=blue]{hyperref}

\usepackage{amssymb,amsmath, amscd}
\usepackage{times, verbatim}
\usepackage{graphicx}
\usepackage[english]{babel}
 \usepackage[usenames, dvipsnames]{color}
\usepackage{amsmath,amssymb,amsfonts}
\usepackage{enumerate}
\usepackage{parskip}
\usepackage[T1]{fontenc}
\usepackage{tikz-cd}
\usetikzlibrary{calc, arrows, decorations.markings} \tikzset{>=latex}
\usepackage{amsmath,amssymb,amsfonts}
\usepackage{amscd, amssymb, latexsym, amsmath, amscd}
\usepackage{csquotes}
\usepackage[all]{xy}
\usepackage{pb-diagram}
\usepackage{enumerate}
\usepackage{ulem}
\usepackage{anysize}
\marginsize{3cm}{3cm}{3cm}{3cm}
\usepackage[colorlinks=true, pdfstartview=FitV, linkcolor=blue, citecolor=blue, urlcolor=blue]{hyperref}
\usepackage{mathdots}
\usepackage{amssymb,amsmath, amscd}
\usepackage{times, verbatim}
\usepackage{graphicx}
\usepackage[english]{babel}
 \usepackage[usenames, dvipsnames]{color}
\usepackage{amsmath,amssymb,amsfonts}
\usepackage{mathtools}
\usepackage{enumerate}
\usepackage{cleveref}
\usepackage[vcentermath]{youngtab}
\usepackage{mathrsfs}
\usepackage{tikz}
\usetikzlibrary{calc, arrows, decorations.markings} \tikzset{>=latex}
\usepackage{amsmath,amssymb,amsfonts}
\usepackage{amscd, amssymb, latexsym, amsmath, amscd}
\usepackage[all]{xy}
\usepackage{pb-diagram}
\usepackage{enumerate}
\usepackage{ulem}
\usepackage{anysize}
\marginsize{3cm}{3cm}{3cm}{3cm}
\input xy
\xyoption{all}
\usepackage{pb-diagram}
\usepackage[all]{xy}
\usepackage{tabularx}
\usepackage{ytableau}
\usepackage{xcolor}

\input xy
\xyoption{all}

\theoremstyle{plain}
\newtheorem{theorem}{Theorem}[section]
\newtheorem{lemma}{Lemma}

\newtheorem{corollary}{Corollary}

\newtheorem{proposition}{Proposition}

\theoremstyle{definition} \theoremstyle{definition}
\newtheorem{remark}{Remark}

\theoremstyle{remark}

\newcommand{\GL}{{\rm GL}}
\newcommand{\C}{\mathbb{C}}

\newcommand{\T}{\underline{t}}
\newcommand{\Sp}{{\rm Sp}}
\newcommand{\SO}{{\rm SO}}

\title{Character values at elements of order 2}

\begin{document}
\author{Chayan Karmakar}

\address{Indian Institute of Technology Bombay, Powai, Mumbai-400076, INDIA}
\email{karmakar.baban28@gmail.com}


\maketitle
    {\hfill \today}
\begin{abstract}
In this paper we compute the character values of highest weight representations for classical groups of types \( A_n \), \( B_n \), \( C_n \), \( D_n \) and the Exceptional group $G_2$ at all conjugacy classes of order 2. We prove that these character values, if nonzero, can be expressed either as a product involving the dimensions of two highest weight representations from classical subgroups, along with an additional constant term, or as an alternating sum of products of the dimensions of two highest weight representations from the classical subgroups, also accompanied by an extra constant term.         
\end{abstract}

\section{Introduction}
In the work \cite{DP} of Dipendra Prasad, he discovered a certain factorization theorem for characters of ${\rm GL}(mn,\C)$ at certain special elements of the diagonal torus, those of the form \[
  \T \cdot c_n :=
  \begin{pmatrix}
    \T & & \\
     & \omega_n\T & \\
    & & \ddots \\
    & & & \omega_n^{n-1}\T
  \end{pmatrix},
\] where $\T=(t_1,t_2,\cdots,t_m) \in (\C^{*})^{m}$ and $\omega_n$ is a primitive $n$th root of unity. D. Prasad proved that the character of a finite dimensional highest weight representation $\pi_{\underline{\lambda}}$ of ${\rm GL}(mn,\C)$ of highest weight $\underline{\lambda}$ at such elements $\T \cdot c_n$ is the product of  characters of certain highest weight representations of ${\rm GL}(m,\C)$ at the element $\T^{n}=(t_1^{n},t_2^{n},\cdots,t_m^{n})\in (\C^{*})^{m}$.  This work of D. Prasad was recently generalized to all classical groups in \cite{AK}.

 A very special case of the works \cite{DP} and \cite{AK} calculates the character of highest weight representations of classical groups at conjugacy classes of elements of order $2$ of the form $(\underbrace{1,\cdots,1}_{k \ {\rm times}},\underbrace{-1,\cdots,-1}_{k\ {\rm times}})$ for ${\rm GL}(2k,\C)$, and for elements of the form \\ $(\underbrace{1,\cdots,1}_{k \ {\rm times}  },\underbrace{-1,\cdots,-1}_{k \ {\rm times} },\underbrace{-1,\cdots,-1}_{ k \ {\rm times} },\underbrace{1,\cdots,1}_{k \ {\rm times}})$ for both the groups ${\rm Sp}(2k,\C)$ and ${\rm SO}(2k,\C)$, along with the element $(\underbrace{1,\cdots,1}_{k \ { \rm times}  },\underbrace{-1,\cdots,-1}_{ k \ { \rm times} },1,\underbrace{-1,\cdots,-1}_{k \ {\rm times} },\underbrace{1,\cdots,1}_{k \ {\rm times} })$ of ${\rm SO}(2k+1,\C)$.\
  
  In this work, we will calculate the character $\Theta_{\lambda}(x_0)$ at any element $x_{0}$ of order $2$ in $G$ when $G$ is a classical group or $G=G_{2}$, generalizing the work of \cite{DP}, \cite{AK} and \cite{N}. 
\section{Some elementay Results for square matices}
In this section we discuss an elementary lemma regarding the determinant of a square matrix generalizing the well-known expression of the determinant of an $n \times n$ matrix in terms of the first column.

We will omit the proof of the following lemma which is a direct consequence of expressing the determinant of a matrix in terms of the linear transformation on the highest exterior power of the underlying vector space.
\begin{lemma}\label{lem 1}
	Let \( B = [b_{ij}] \) be an \( n \times n \) matrix, and let \( S=\{ s_1 < \cdots <s_k\} \) be a  subset of $\{1,2,\cdots,n \}$ of $k$ elements. Let $C_{S}$ be the set of all $k \times k$ submatrices of $B$ taken from the columns indexed by $S$. For each $k \times k$ matrix $X \in C_{S}$, let $X^{'}$ be the $(n-k) \times (n-k)$ matrix consisting of complementary rows and columns, Then:
	
	\[
	\det(B)= \sum_{X \in C_{S}} \epsilon_{X} \det(X) \det(X^{'}),
	\]
	where $\epsilon_{X}$ is the sign of the \enquote{shuffle} permutation: if $X$ corresponds to rows $\{ r_{1} < \cdots < r_{k} \}$  and $X^{'}$ corresponds to rows $\{t_{1}<\cdots < t_{n-k}\}$, then 
	\[
	\epsilon_{X}(e_1 \wedge \cdots \wedge e_{n})=e_{r_1} \wedge \cdots e_{r_{k}} \wedge e_{t_{1}} \wedge \cdots \wedge e_{t_{n-k}}.
	\]   
	
\end{lemma}

The following corollaries are special case of Lemma \ref{lem 1}.

\begin{corollary}\label{cor 1}
Let $A$ be a square matrix of order $n$ and $B$ is a submatrix of size $a \times b$ with $a+b>n$, and B=0, then $ \det A =0$.	
\end{corollary}

\begin{corollary}\label{cor 2}
Let $A$ be a square matrix of order $n$ and $B$ is a submatrix of size $a \times b$ with $a+b=n$, and $B=0$, then we have 
$$\det A=\det A_1 \times \det A_2,$$
where $A_{1}$  and $A_{2}$ are  complementary square submatrices of order $a$ and $b$ respectively.	
\end{corollary}
 
\section{Definitions and Notations}\label{sec: defn}
We begin with some general notation to be used throughout the paper.

Let \( G \) be one of the classical groups: \(\GL(n, \mathbb{C})\), \(\Sp(2n, \mathbb{C})\), \(\SO(2n, \mathbb{C})\), or \(\SO(2n+1, \mathbb{C})\). We consider a highest weight representation of \( G \) characterized by the highest weight \(\underline{\lambda}= (\lambda_1, \lambda_2, \ldots, \lambda_n)\), where \( \lambda_1 \geq \lambda_2 \geq \cdots \geq \lambda_n \geq 0 \).

For  $0 \leq i \leq 1$, we define the set \(\eta_i(\underline{\lambda})\) as follows:
\[
\eta_i(\underline{\lambda}) = \{ a \in \underline{\lambda} + \rho \mid a \equiv i \mod 2 \},
\]
where $a \in \underline{\lambda} + \rho$ means, $a$ appears as some coordinate of $\underline{\lambda}+\rho$, \(\rho\) is the half the sum of the positive roots of \( G \) and $\rho$ is given by:

\begin{align*}
\rho=\begin{cases}
\rho_{n}=(n-1,n-2,\cdots,0) & \text{if} \ G=\GL(n,\C), \\
\rho^{n}=(n,n-1,\cdots,1) & \text{if} \  G=\Sp(2n,\C), \\
\overline{\rho}_{2n}=(n-1,n-2,\cdots,0) & \text{if} \  G=\SO(2n,\C), \\
\overline{\rho}_{2n+1}=(n-1/2,n-3/2,\cdots,1/2) & \text{if} \  G=\SO(2n+1,\C). \\
	\end{cases} \label{0.5}
\end{align*}

When \(0 \leq i \leq 1\), we define \(\eta_i(\underline{\lambda})/2 = \left\{ \frac{a}{2} \mid a \in \eta_i(\underline{\lambda}) \right\}\) and \([\eta_1(\underline{\lambda}) - 1]/2 = \left\{ \frac{b - 1}{2} \mid b \in \eta_1(\underline{\lambda}) \right\}\).   

For any positive integer $n$, let us define $[n]=\{1,2,\cdots,n\}$. For any $n$-tuple $a=(a_1,a_2,\cdots,a_{n}) \in (\C^{\ast})^{n}$, define $a^{t}=(a_{n},a_{n-1},\cdots,a_{1})$.

\section{The General Linear Group ${\rm GL}(n,\C)$}
\begin{theorem}\label{thm 1}
Let $G={\rm GL}(n,\C)$. For each integer $0<k \leq n/2$, consider the diagonal element $$C_{n-k,k}=(\underbrace{1,1,\cdots,1}_{n-k \ {\rm times}},\underbrace{-1,-1,\cdots,-1}_{k \ {\rm times}}).$$ These are elements of order 2 in $G$. Let $S_{(\underline{\lambda})}\C^{n}$ be the highest weight representation of ${\rm GL}(n,\C)$ with highest weight $\underline{\lambda}=(\lambda_1,\lambda_2,\cdots,\lambda_n)$ with $\lambda_1 \geq \lambda_2 \geq \cdots \geq \lambda_{n}$. Assume after possibly multiplying $S_{(\underline{\lambda})}\C^{n}$ by the determinant character that $\# \eta_{0}(\underline{\lambda}) \geq \# \eta_{1}(\underline{\lambda})$. Then :
\begin{itemize}
    \item [A.] The character $\Theta_{\underline{\lambda}}(C_{n-k,k})=0$ if  $\#\eta_{0}(\underline{\lambda})>n-k$. \\
    \item [B.] For $\#\eta_{0}(\underline{\lambda})=n-k$, we have the following factorization:
\[
\Theta_{\underline{\lambda}}(C_{n-k,k})= \pm 2^{c(k)} \cdot \dim(S_{(\lambda^{0})}\C^{n-k} \otimes S_{(\lambda^{1})}\mathbb{C}^{k}),
\]
where $c(k)=\binom{n-2k}{2}$ and the highest weights $\lambda^{0}$ and $\lambda^{1}$  are given by: \\
$
\begin{cases}
	\lambda^{0}+\rho_{n-k}=\eta_{0}(\underline{\lambda})/2, \\
	\lambda^{1}+\rho_{k}=[\eta_{1}(\underline{\lambda})-1]/2,
\end{cases}
$ \\ \\ 
 where $\rho_{k}=(k-1,k-2,\cdots,0).$ 
 By convention, the elements of $\eta_{i}(\underline{\lambda})$ are arranged in the decreasing order.  \\ 

\item[C.] Let $\underline{\lambda}+\rho_{n}=(l_{1},l_{2},\cdots,l_{n})$. If $\# \eta_{0}(\underline{\lambda})<n-k$, we have the following alternating sum:
\[
\Theta_{\underline{\lambda}}(C_{n-k,k})=\pm \frac{ \sum_{S \in I}\epsilon_{S}\cdot \rm dim (S_{(\lambda_{S})}\C^{k} \otimes S_{(\lambda_{S^{'}})}\C^{n-k}) }{2^{k(n-k-1)}},
\]
where $S^{'}=[n] \backslash S$, $I=\{ \{i_1,i_2,\cdots,i_k\} \subset [n] \  \vert  \   l_{i_{j}} \in \eta_{1}(\underline{\lambda}) \  \forall \ 1 \leq j \leq k \}$, $\epsilon_{S}$ is the sign of the \enquote{shuffle} permutation corresponding to the sets of rows indexed by $S$ and $S^{'}$.
The highest weights $\lambda_{S}$ of $\GL(k,\C)$ and $\lambda_{S^{'}}$ of $\GL(n-k,\C)$ are given by:
\small
\begin{enumerate}
	\item[(i.)]$\lambda_{S}+\rho_k=(l_{s} \vert s \in S)$. 
	
	\item[(ii.)]$\lambda_{S^{'}}+\rho_{n-k}=(l_{s} \vert s \in S^{'})$.
\end{enumerate}
\noindent There are a total of ${\# \eta_{1}(\underline{\lambda}) \choose k}$ terms in the above summation.
\end{itemize}

 \begin{proof}
The proof of this theorem will be a direct application of the Weyl character formula. Note that the elements $C_{n-k,k}$ are singular elements of ${ \rm GL}(n,\C)$, hence the Weyl denominator is zero. Therefore $\Theta_{\lambda}(C_{n-k,k})$ is calculated by taking a limit $\Theta_{\lambda}(C_{n-k,k}(\epsilon))$, where $C_{n-k,k}(\epsilon)$ are elements of ${\rm GL}(n,\C)$ which converge to $C_{n-k,k}$ as $\epsilon \rightarrow 0$. There are many options to consider $C_{n-k,k}(\epsilon)$, but we take  $$C_{n-k,k}(\epsilon)=(x_{n-k}(\epsilon),x_{n-k+1}(\epsilon),\cdots,x_{1}(\epsilon),-x_{k}(\epsilon),-x_{k-1}(\epsilon),\cdots,-x_{1}(\epsilon)),$$ where  \[ x_{j}(\epsilon)=1+j\epsilon , 1 \leq j \leq n-k. \] These are regular elements for each $\epsilon >0$ and therefore the Weyl denominator is nonzero. We will use the Weyl character formula as the quotient of two determinants given by
$$
\Theta_{\underline{\lambda}}(x_1,x_2,\cdots,x_n)=\frac{\det (x_j^{\lambda_{i}+n-i} )}{\Delta(x_1,x_2,\cdots,x_n)},
$$ where $\Delta$ denote the Vandermonde determinant.  
We have $$\lim_{\epsilon \rightarrow 0}\frac{A_{\underline{\lambda}+\rho_{n}}(C_{n-k,k}(\epsilon))}{A_{\rho_{n}}(C_{n-k,k}(\epsilon))}=\Theta_{\underline{\lambda}}(C_{n-k,k}).$$ 

So we will begin by calculating the Weyl denominator $A_{\rho_{n}}(C_{n-k,k}(\epsilon))$. \\ 
\begin{center}
\textbf{Calculation of the Weyl denominator at $C_{n-k,k}(\epsilon)$}	
\end{center}

The Weyl denominator at $C_{n-k,k}(\epsilon)$ is given by:
\begin{align}
	\small
	\nonumber
	A_{\rho_n}(C_{n-k,k}(\epsilon))&=\Delta(C_{n-k,k}(\epsilon)), \label{2.0} \\
	&=\Delta(x_{n-k}(\epsilon),x_{n-k-1}(\epsilon),\cdots,x_{1}(\epsilon),-x_k(\epsilon),-x_{k-1}(\epsilon),\cdots,-x_{1}(\epsilon)), \\
	&=(-1)^{\frac{k(k-1)}{2}}\Delta(C(\epsilon)) \times \Delta(D(\epsilon)) \times \prod^{n-k}_{s=1}\prod^{k}_{t=1}(x_{s}(\epsilon)+x_{t}(\epsilon)), \label{2.1}\\
	&=(-1)^{\frac{k(k-1)}{2}}\frac{\Delta(C^{2}(\epsilon)) \times \Delta(D^{2}(\epsilon)) \times  \prod^{k}_{s=1}2x_s(\epsilon)}{\prod_{\{ s > t \vert k+1 \leq s,t \leq n-k \}}(x_s(\epsilon)+x_{t}(\epsilon))}\label{2.2} , \\
	&=(-1)^{\frac{k(k-1)}{2}}\frac{2^{k}\Delta(C^{2}(\epsilon)) \times \Delta(D^{2}(\epsilon)) \times \prod^{k}_{s=1}x_{s}(\epsilon)}{\prod_{\{ s > t \vert k+1 \leq s,t \leq n-k \}}(x_{s}(\epsilon)+x_{t}(\epsilon))}\label{2.3},
\end{align}
where  $$C(\epsilon)=(x_{n-k}(\epsilon),x_{n-k-1}(\epsilon),\cdots,x_{1}(\epsilon)) \in {\rm GL}(n-k,\C)$$ and $$D(\epsilon)=(x_{k}(\epsilon),x_{k-1}(\epsilon),\cdots,x_{1}(\epsilon)) \in  {\rm GL}(k,\C);$$ and where $C^{2}(\epsilon)$ and $D^{2}(\epsilon)$ are the squares of the corresponding matrices.

Now to show that $\Theta_{\underline{\lambda}}(C_{n-k,k})=0$, it is sufficient to consider the conditions under which the Weyl numerator $A_{\underline{\lambda}+\rho_{n}}(C_{n-k,k}(\epsilon))$ is identically zero which is what we do now.  

Let $\underline{\lambda}+\rho_n=(l_1,l_2,\cdots,l_n)$ with $l_i=\lambda_i+n-i$. 
Let $\#\eta_{0}(\underline{\lambda})=n-k+s$ with $ 1 \leq s \leq k$. Interchanging the rows of \(A_{\underline{\lambda}+\rho_{n}}(C_{n-k,k}(\epsilon))\) changes the determinant only by a sign. Therefore, without loss of generality, we can assume that  \(\eta_{0}(\underline{\lambda})=\{l_1, l_2, \ldots, l_{n-k+s}\}\).  
Now the Weyl numerator at $C_{n-k,k}(\epsilon)$ is given by 
 \[
 A_{\underline{\lambda}+\rho_{n}}(C_{n-k,k}(\epsilon))= \det \begin{pmatrix}
 	x_{n-k}(\epsilon)^{l_1} & \cdots & x_{1}(\epsilon)^{l_1} & x_{k}(\epsilon)^{l_1} & \cdots & x_{1}(\epsilon)^{l_1} \\ \\
 	x_{n-k}(\epsilon)^{l_2} & \cdots & x_{1}(\epsilon)^{l_2} & x_{k}(\epsilon)^{l_2} & \cdots & x_{1}(\epsilon)^{l_2} \\
 	\vdots & \ddots & \vdots & \vdots & \ddots & \vdots \\
 	x_{n-k}(\epsilon)^{l_{n-k+s}} & \cdots & x_{1}(\epsilon)^{l_{n-k+s}} & x_{k}(\epsilon)^{l_{n-k+s}} & \cdots & x_{1}(\epsilon)^{l_{n-k+s}} \\
 	x_{n-k}(\epsilon)^{l^{'}_{s+1}} & \cdots & x_{1}(\epsilon)^{l^{'}_{s+1}} & -x_{k}(\epsilon)^{l^{'}_{s+1}} & \cdots & -x_{1}(\epsilon)^{l^{'}_{s+1}} \\
 	\vdots & \ddots & \vdots & \vdots & \ddots & \vdots \\
 	x_{n-k}(\epsilon)^{l^{'}_{n}} & \cdots & x_{1}(\epsilon)^{l^{'}_{n}} & -x_{k}(\epsilon)^{l^{'}_{n}} & \cdots & -x_{1}(\epsilon)^{l^{'}_{n}} \\
 \end{pmatrix},
 \]
 where for $1 \leq i \leq k-s$ , we define $l_{s+i}^{'}=l_{n-k+s+i}$. For each $1 \leq i \leq k$, we will apply the following set of successive column operations on the Weyl numerator $A_{\underline{\lambda}+\rho_{n}}(C_{n-k,k}(\epsilon))$:
 \begin{itemize}
 	\item[(1)]\( C_{n-k+i} \rightarrow C_{n-k+i}-C_{n-2k+i}, \)
 	\item[(2)]\(C_{n-k+i} \rightarrow -\frac{1}{2}C_{n-k+i}. \)  
 \end{itemize}
 By applying these successive operations, the matrix corresponding to the Weyl numerator $A_{\underline{\lambda}+\rho_{n}}(C_{n-k,k}(\epsilon))$ transforms into the following matrix of the form:
 $$
 (-2)^{k} \times \begin{pmatrix}
 	A_{(n-k+s) \times (n-k)}(\epsilon) & \textbf{0}_{(n-k+s) \times k}\\
 	C_{(k-s) \times (n-k)}(\epsilon) & D_{(k-s) \times k}(\epsilon)
 \end{pmatrix},
 $$
 where  $$A_{(n-k+s) \times (n-k)}(\epsilon)=\begin{pmatrix}
 	x_{n-k}(\epsilon)^{l_1} & \cdots & x_{1}(\epsilon)^{l_1}  \\ \\
 	x_{n-k}(\epsilon)^{l_2} & \cdots & x_{1}(\epsilon)^{l_2}  \\
 	\vdots & \ddots & \vdots  \\
 	x_{n-k}(\epsilon)^{l_{n-k+s}} & \cdots & x_{1}(\epsilon)^{l_{n-k+s}}  \\
 \end{pmatrix}$$ and
 $$D_{(k-s) \times k}(\epsilon)=\begin{pmatrix}
 	x_{k}(\epsilon)^{l^{'}_{s+1}} & \cdots & x_{1}(\epsilon)^{l^{'}_{s+1}} \\ \\
 	x_{k}(\epsilon)^{l^{'}_{s+2}} & \cdots & x_{1}(\epsilon)^{l^{'}_{s+2}} \\
 	\vdots & \ddots & \vdots \\
 	x_{k}(\epsilon)^{l^{'}_{n}} & \cdots & x_{1}(\epsilon)^{l^{'}_{s+1}} \\
 \end{pmatrix}.$$ Note that for the zero matrix $\textbf{0}_{(n-k+2s) \times (k-s)}$, we have $(n-k+s)+k=n+s>n$. Therefore by Corollary \ref{cor 1} we deduce that  $A_{\underline{\lambda}+\rho_{n}}(C_{n-k,k}(\epsilon))=0$.   
 Thus $A_{\underline{\lambda}+\rho_{n}}(C_{n-k,k})=0$ whenever $\#\eta_{0}(\underline{\lambda})>n-k$. 
 
This completes the proof of Part (A) of Theorem \ref{thm 1}.

Next we turn to the proof of Part$(B)$ of the theorem, thus we assume that \( \#\eta_{0}(\underline{\lambda}) = n-k \).  
 
Without loss of generality, we can assume that \( \eta_{0}(\underline{\lambda})=\{ l_1, l_2, \ldots, l_{n-k} \}. \) 
The Weyl numerator at $C_{n-k,k}(\epsilon)$ is given by
\[
\det
\begin{pmatrix}
	x_{n-k}(\epsilon)^{l_1} & \cdots & x_{1}(\epsilon)^{l_1} & x_{k}(\epsilon)^{l_1} & \cdots & x_{1}(\epsilon)^{l_1} \\ \\
	x_{n-k}(\epsilon)^{l_2} & \cdots & x_{1}(\epsilon)^{l_2} & x_{k}(\epsilon)^{l_2} & \cdots & x_{1}(\epsilon)^{l_2} \\
	\vdots & \ddots & \vdots & \vdots & \ddots & \vdots \\
	x_{n-k}(\epsilon)^{l_{n-k}} & \cdots & x_{1}(\epsilon)^{l_{n-k}} & x_{k}(\epsilon)^{l_{n-k}} & \cdots & x_{1}(\epsilon)^{l_{n-k}} \\ \\
	x_{n-k}(\epsilon)^{l_{n-k+1}} & \cdots & x_{1}(\epsilon)^{l_{n-k+1}} & -x_{k}(\epsilon)^{l_{n-k+1}} & \cdots & -x_{1}(\epsilon)^{l_{n-k+1}} \\
	\vdots & \ddots & \vdots & \vdots & \ddots & \vdots \\
	x_{n-k}(\epsilon)^{l_{n}} & \cdots & x_{1}(\epsilon)^{l_{n}} & -x_{k}(\epsilon)^{l_{n}} & \cdots & -x_{1}(\epsilon)^{l_{n}} \\
\end{pmatrix},
\]
For each $1 \leq i \leq k$, we will apply the following set of column operations on $A_{\underline{\lambda}+\rho_n}(C_{n-k,k}(\epsilon))$:
\begin{itemize}
	\item[(i)]\( C_{n-k+i} \rightarrow C_{n-k+i}-C_{n-2k+i}, \)
	\item[(ii)]\(C_{n-k+i} \rightarrow -\frac{1}{2}C_{n-k+i}. \)
\end{itemize} Applying these operations $A_{\underline{\lambda}+\rho_{n}}(C_{n-k,k}(\epsilon))$ can be transformed into the block diagonal form given by 
\[
(-2)^{k} \times \det
\begin{pmatrix}
	E_{(n-k) \times (n-k)}(\epsilon) & 0_{(n-k) \times k} \\
	K_{k \times (n-k)}(\epsilon) & F_{k \times k}(\epsilon)
\end{pmatrix},
\]
where $E_{(n-k) \times (n-k)}(\epsilon)=(x_{n-k-j+1}(\epsilon)^{l_i})_{ij}$ and 
$F_{k \times k}(\epsilon)=(x_{k-j+1}(\epsilon)^{(l_{n-k+i})})_{ij}$.
Let us consider the matrix \[\begin{pmatrix}
	E_{(n-k) \times (n-k)}(\epsilon) & 0_{(n-k) \times k} \\
	K_{k \times (n-k)}(\epsilon) & F_{k \times k}(\epsilon)
\end{pmatrix}. \]
It has the zero submatrix $0_{(n-k) \times k}$ and also $(n-k)+k=n$. Therefore by Corollary \ref{cor 2} we have
\begin{align}
	\nonumber
	A_{\underline{\lambda}+\rho_n}(C_{n-k,k}(\epsilon))&=\pm 2^{k} \det (E_{n-k \times n-k}(\epsilon)) \times \det (F_{k \times k}(\epsilon)), \\
	&=\pm 2^{k} \det (E_{n-k \times n-k}(\epsilon)) \times \det (F^{'}_{k \times k}(\epsilon)) \times \prod_{j=1}^{k}x_{j}(\epsilon),\label{6}
\end{align}
where $E_{(n-k) \times (n-k)}(\epsilon)$ is defined earlier and  
$F^{'}_{k \times k}(\epsilon)=(x_{k-j+1}(\epsilon)^{(l_{n-k+i}-1)})_{ij}$.

Note that we have $l_{m}\equiv\begin{cases}
	0 \mod 2 & \text{if} \ m=1, \cdots,n-k \\
	1 \mod 2 & \text{if} \ m=n-k+1,\cdots,n
\end{cases}$ \\ 

Now from equations \eqref{2.3}, \eqref{6} and by the applications of Weyl character formula we have:

\begin{align}
\nonumber
\Theta_{\underline{\lambda}}(C_{n-k,k}(\epsilon))&=\pm \prod_{\{ s > t \vert k+1 \leq s,t \leq n-k \}}(x_{s}(\epsilon)+x_{t}(\epsilon)) \times \frac{\det (E_{n-k \times n-k}(\epsilon))}{\Delta(C^{2}(\epsilon))} \times \frac{\det (F^{'}_{k \times k}(\epsilon))}{(D^{2}(\epsilon))}, \\
&=\pm \prod_{\{ s > t \vert k+1 \leq s,t \leq n-k \}}(x_{s}(\epsilon)+x_{t}(\epsilon)) \times \Theta_{\lambda^{0}}(C^{2}(\epsilon)) \times \Theta_{\lambda^{1}}(D^{2}(\epsilon)),\label{7}
\end{align}  

where the highest weights $\lambda^{0}$ of $\GL(n-k,\C)$ and  $\lambda^{1}$ of $\GL(k,\C)$ are defined by :

\begin{itemize}
	\item [(i.)] $\lambda^{0}+\rho_{n-k}=\eta_{0}(\lambda)/2=(\frac{l_1}{2},\frac{l_2}{2},\cdots,\frac{l_{n-k}}{2})$.
	
	\item [(ii.)] $\lambda^{1}+\rho_{k}=[\eta_{1}(\lambda)-1]/2=(\frac{l_{n-k+1}-1}{2},\frac{l_{n-k+2}-1}{2},\cdots,\frac{l_{n}-1}{2}).$
	
\end{itemize}

Now from equation \eqref{7} we get 
\begin{align}
\nonumber
\Theta_{\underline{\lambda}}(C_{n-k,k})&=\lim_{\epsilon \rightarrow 0}\Theta_{\underline{\lambda}}(C_{n-k,k}(\epsilon)), \\
&=\prod_{\{ i > j \vert k+1 \leq i,j \leq n-k \}}\lim_{\epsilon \rightarrow 0}(x_{s}(\epsilon)+x_{t}(\epsilon)) \times \lim_{\epsilon \rightarrow 0} \Theta_{\lambda^{0}}(C^{2}(\epsilon)) \times \lim_{\epsilon \rightarrow 0} \Theta_{\lambda^{1}}(D^{2}(\epsilon)), \\
&=2^{{n-2k \choose 2}} \times {\rm dim}(S_{(\lambda^{0})}\C^{n-k}) \times {\rm dim}(S_{(\lambda^{1})}\C^{k}), \\
&=2^{c(k)} \times {\rm dim}( S_{(\lambda^{0})}\C^{n-k} \otimes S_{(\lambda^{1})}\C^{k}),
\end{align}
where $c(k)={n-2k \choose 2}$. 
This completes the proof of Part (B) of Theorem \ref{thm 1}. \\ \\
Next we will prove part (C) of this theorem. Let  $\#\eta_{0}(\underline{\lambda})=m$ with $k+1 \leq m \leq n-k-1$. Without loss of generality let us assume that $\eta_0(\underline{\lambda})=\{l_1,l_2,\cdots,l_{m}\}$. For each $1 \leq i \leq k$, we will apply the following set of successive column operations on the Weyl numerator $A_{\lambda+\rho_{n}}(C_{n-k,k}(\epsilon))$:
\begin{itemize}
	\item[(1)]\( C_{n-k+i} \rightarrow C_{n-k+i}-C_{n-2k+i}, \)
	\item[(2)]\(C_{n-k+i} \rightarrow -\frac{1}{2}C_{n-k+i}, \)
\end{itemize}
By applying these successive operations the Weyl numerator $A_{\underline{\lambda}+\rho_{n}}(C_{n-k,k}(\epsilon))$ transforms into 
$$
(-2)^{k} \times \det(B(\epsilon)),
$$
where  

 $$B(\epsilon)=\det \begin{pmatrix}
	x_{n-k}(\epsilon)^{l_1} & \cdots & x_{1}(\epsilon)^{l_1} &0 & \cdots & 0 \\ \\
	x_{n-k}(\epsilon)^{l_2} & \cdots & x_{1}(\epsilon)^{l_2} &0 & \cdots & 0 \\
	\vdots & \ddots & \vdots & \vdots & \ddots & \vdots \\
	x_{n-k}(\epsilon)^{l_m} & \cdots & x_{1}(\epsilon)^{l_m} &0 & \cdots & 0 \\ \\
	x_{n-k}(\epsilon)^{l_{m+1}} & \cdots & x_{1}(\epsilon)^{l_{m+1}} & x_k(\epsilon)^{l_{m+1}} & \cdots & x_1(\epsilon)^{l_{m+1}} \\
	\vdots & \ddots & \vdots & \vdots & \ddots & \vdots \\
	x_{n-k}(\epsilon)^{l_{n}} & \cdots & x_{1}(\epsilon)^{l_{n}} & x_k(\epsilon)^{l_{n}} & \cdots & x_1(\epsilon)^{l_{n}} \\
\end{pmatrix}.$$

Applying Lemma \ref{lem 1} on $B(\epsilon)$ by choosing the subset $H=\{n-k+1,\cdots,n\}$, which corresponds to the last $k$ columns, we get 

\begin{align}
	B(\epsilon) &= \sum_{\substack{ S \subset H_m \\ \\ \# S = k }} \epsilon_{S} \cdot \det(M_{S,H}(\epsilon)) \det(M_{S',H'}(\epsilon)),
\end{align}
where \( H_m = \{ m+1, m+2, \dots, n \} \), 
$\epsilon_{S}$ is the sign of the \enquote{shuffle} permutation corresponding to the sets of rows indexed by $S$ and $S^{'}$.
\( M_{S,H}(\epsilon) \) denotes the \( k \times k \) submatrix of the matrix \( B(\epsilon) \) obtained by choosing rows with indices in \( S \) and columns with indices in \( H \), and \( M_{S',H'}(\epsilon) \) denotes the submatrix in \( B(\epsilon) \) having rows with indices in $S^{'}=[n] \setminus S$ and columns with indices in $H^{'}=[n] \setminus H$. In particular $$ M_{S,H}(\epsilon)=(x^{l_{u}}_{j}(\epsilon))_{u \in S, j \in [k]}$$ and $$M_{S',H'}(\epsilon)=(x^{l_{u}}_{j}(\epsilon))_{u \in S', j \in [n-k]},$$ where for any integer $w$, we have defined $[w]$ in Section \ref{sec: defn}.

Now the Weyl numerator is given by
\begin{align}
\nonumber
A_{\underline{\lambda}+\rho_{n}}(C_{n-k,k}(\epsilon))&=(-2)^{k} \times B(\epsilon), \\
&=(-2)^{k} \times \sum_{\substack{ S \subset H_m \\ \\ \# S = k }}
\epsilon_{S} \cdot \text{det}(M_{S,H}(\epsilon))\text{det}(M_{S^{'},H^{'}}(\epsilon)).\label{14.1}
\end{align}

Now from equations \eqref{2.1}, \eqref{14.1} and by the applications of Weyl character formula we get
\small
\begin{align}
\Theta_{\underline{\lambda}}(C_{n-k,k}(\epsilon))&=\pm2^{k} \times  \sum_{\substack{ S \subset H_m \\ \\ \# S = k }}
\frac{\epsilon_{S}}{\prod^{n-k}_{s=1}\prod^{k}_{t=1}(x_s(\epsilon)+x_t(\epsilon))} \times   \frac{\det(M_{S,H}(\epsilon))}{\Delta(C(\epsilon))} \times  \frac{\det(M_{S^{'},H}(\epsilon))}{\Delta(D(\epsilon))},\label{12} \\
&=\pm 2^{k} \times \sum_{\substack{ S \subset H_m \\ \\ \# S = k }}\frac{\epsilon_{S}}{\prod^{n-k}_{s=1}\prod^{k}_{t=1}(x_s(\epsilon)+x_t(\epsilon))} \times  \Theta_{\lambda_{S}}(C(\epsilon)) \times \Theta_{\lambda_{S^{'}}}(D(\epsilon)),\label{14}
\end{align}

where the highest weights $\lambda_{S}$ of $\GL(k,\C)$ and $\lambda_{S^{'}}$ of $\GL(n-k,\C)$ are defined by:

\begin{itemize}
	\item [1.] $\lambda_{S}+\rho_{k}=(l_{s} \vert u \in S)$.
	
	\item [2.] $\lambda_{S^{'}}+\rho_{n-k}=(l_{u} \vert u \in S^{'})$. 
\end{itemize} 

From equation \eqref{14} we have 
\small
\begin{align*}
\Theta_{\underline{\lambda}}(C_{n-k,k})&=\lim_{\epsilon \rightarrow 0}\Theta_{\underline{\lambda}}(C_{n-k,k}(\epsilon)), \\
&=\pm \sum_{\substack{ S \subset H_m \\ \\ \# S = k }}\frac{2^{k}\epsilon_{S}}{\prod^{n-k}_{s=1}\prod^{k}_{t=1}\lim_{\epsilon \rightarrow 0}(x_s(\epsilon)+x_t(\epsilon))} \times  \lim_{\epsilon \rightarrow 0}\Theta_{\lambda_{S}}(C(\epsilon)) \times \lim_{\epsilon \rightarrow 0}\Theta_{\lambda_{S^{'}}}(D(\epsilon)), \\
&=\pm \sum_{\substack{ S \subset H_m \\ \\ \# S = k }}2^{k} \epsilon_{S} \cdot \frac{{\rm dim}(S_{(\lambda_{S})}\C^{k}){\rm dim}(S_{(\lambda_{S^{'}})}\C^{n-k})}{2^{k(n-k)}}, \\
&=\pm \frac{1}{2^{k(n-k-1)}} \sum_{S \in I}\epsilon_{S} \cdot {\rm dim}(S_{(\lambda_{S})}\C^{k} \otimes S_{(\lambda_{S^{'}})}\C^{n-k}),
\end{align*}

where $I=\{S \subset H_{m} \vert \#S=k \}$ and $H_{m}=\{j \mid l_{j} \in \eta_{1}(\underline{\lambda}) \}$. \\ \\
This completes the proof of Part (C) of Theorem \ref{thm 1}.
\end{proof} 
\end{theorem}
\begin{remark}
Observe that $c(k)$ is 0 if we have $n=2k$ and $2k+1$ respectively. Therefore when $n=2k$ or $n=2k+1$, the character value at $C_{n-k,k}$ can be written up to a sign as the dimension of an highest weight representation of ${\rm GL}(n-k,\C) \times {\rm GL}(k,\C)$. Thus Theorem \ref{thm 1} generalizes Theorem 2 of \cite{DP} and Theorem 2.5 of \cite{AK} in the particular case of elements of order $2$. 
\end{remark} 

 \section{The Symplectic Group ${\rm Sp}(2n,\C)$} \label{Sec:Symplectic}
For any $n$-tuple $a=(a_1,a_2,\cdots,a_{n}) \in (\C^{\ast})^{n}$, let us recall the notation  $a^{t}$ defined in Section \ref{sec: defn}.
\begin{theorem}\label{thm 2}
Let $G={\rm Sp}(2n,\C)$. For each integer $0<k \leq n/2$, consider $D_{n-k,k}=(C_{n-k,k},(C^{t}_{n-k,k})^{-1})$, where $C_{n-k,k}=(\underbrace{1,1,\cdots,1}_{n-k \ {\rm times}},\underbrace{-1,-1,\cdots,-1}_{k \ {\rm times}})$. These are elements of order 2. Let $S_{\langle \underline{\lambda} \rangle}\C^{2n}$ be the highest weight representation of ${\rm Sp}(2n,\C)$ with highest weight $\underline{\lambda}=(\lambda_1,\lambda_2,\cdots,\lambda_n)$. Then:
\begin{itemize}
    \item[A.] The character $\Theta_{\underline{\lambda}}(D_{n-k,k})=0$ if either $\#\eta_{0}(\underline{\lambda})>n-k$ or $\#\eta_{1}(\underline{\lambda})>n-k$. \\
    \item[B.] For $\#\eta_{0}(\underline{\lambda})=n-k$ or $\#\eta_{0}(\underline{\lambda})=k$, we have the following factorization:
\[
\noindent \Theta_{\underline{\lambda}}(D_{n-k,k})= 
\begin{cases}
	\pm 2^{d(k)} \cdot \dim(S_{\langle \lambda^{0} \rangle} \mathbb{C}^{2n-2k} \bigotimes S_{[\lambda^{1}]} \mathbb{C}^{2k+1}) , & \text{if } \#\eta_{0}(\lambda) = n-k, \\[1ex]
	\pm 2^{d(k)} \cdot \dim(S_{\langle \lambda^{0} \rangle} \mathbb{C}^{2k} \otimes S_{[\lambda^{1}]} \mathbb{C}^{2n-2k+1}) , & \text{if } \#\eta_{0}(\lambda) = k,	
\end{cases}
\]
where $d(k)={(n-2k)^{2}}$ and the highest weights $\lambda^{0}$ and $\lambda^{1}$ are defined by: \\ \\
$
\begin{cases}
	\begin{cases}
		\lambda^{0}+\rho^{n-k}=\eta_{0}(\underline{\lambda})/2, \\
	\lambda^{1}+\overline{\rho}_{2k+1}=\eta_{1}(\underline{\lambda})/2.
	\end{cases} & \text{if} \ \#\eta_{0}(\underline{\lambda})=n-k, \\ \\
	\begin{cases}
			\lambda^{0}+\rho^{k}=\eta_{0}(\underline{\lambda})/2, \\
		\lambda^{1}+\overline{\rho}_{2n-2k+1}=\eta_{1}(\underline{\lambda})/2.
	\end{cases} & \text{if} \ \#\eta_{0}(\underline{\lambda})=k ,
\end{cases}
$ \\ \\
where $\rho^{n}$ and $\overline{\rho}_{2k+1}$ are the half the sum of the positive roots of ${\rm Sp}(2n,\C)$ and ${\rm SO}(2k+1,\C)$ respectively, for example $\rho^{n}=(n,n-1,\cdots,1)$.

\item[C.] Let $\underline{\lambda}+\rho^{n}=(l_1,l_2,\cdots,l_n)$. For $k< \#\eta_{0}(\lambda)<n-k$, we have the following alternating sum:
\small
\[
\Theta_{\underline{\lambda}}(D_{n-k,k})=\pm \frac{\sum_{S \in I}\epsilon_{S}\cdot {\rm dim}(S_{ \langle \lambda_{S} \rangle}\mathbb{C}^{2n-2k} \otimes S_{ \langle \lambda_{S^{'}} \rangle}\mathbb{C}^{2k})}{2^{k(2n-2k-1)}},
\]
where $S^{'}=[n] \backslash S$, $I=\{ \{i_1,i_2,\cdots,i_k\} \subset [n] \vert  \   l_{i_{j}  } \in \eta_{1}(\lambda)  \  \forall \ 1 \leq j \leq k\}$, $\epsilon_{S}$ is the sign of the \enquote{shuffle} permutation corresponding to the sets of rows indexed by $S$ and $S^{'}$.
The highest weights $\lambda_{S}$ of ${\rm Sp}(2k,\C)$ and $\lambda_{S^{'}}$ of ${\rm Sp}(2n-2k,\C)$ are given by:
\small
\begin{enumerate}
	\item[(i.)]$\lambda_{S}+\rho^{k}=(l_{s} \vert s \in S)$. 
	
	\item[(ii.)]$\lambda_{S^{'}}+\rho^{n-k}=(l_{s} \vert s \in S^{'})$.
\end{enumerate}
\noindent There are a total of ${\# \eta_{1}(\underline{\lambda}) \choose k}$ terms in the above summation.
\end{itemize}
\begin{proof}
The proof of this theorem is analogous to the proof of the corresponding theorem for $\GL(n,\C)$ in which we replace the Weyl numerator (and the Weyl denominator) $\det(x^{\lambda_i+n-i}_{j})$ by $\det(x^{\lambda_i+n+1-i}_{j}-x^{-\lambda_i-n-1+i}_{j})$.
\end{proof}
\end{theorem}
\begin{remark}
	For $n=2k$, we have $d(k)=0$. So that whenever $n=2k$, the character value at $D_{n-k,k}$  can be written up to a sign either as the dimension of a highest weight representation of ${\rm Sp}(2n-2k,\C) \times {\rm SO}(2k+1,\C)$ or as the dimension of a highest weight representation of ${\rm Sp}(2k,\C) \times {\rm SO}(2n-2k+1,\C)$. This shows that Theorem \ref{thm 2} generalizes Theorem 2.11 from \cite{AK} in the particular case of elements of order $2$.     
\end{remark}
\section{The even Orthogonal Group ${\rm SO}(2n,\C)$}
 For $a=(a_1,a_2,\cdots,a_{n}) \in (\C^{\ast})^{n}$, let us recall the notation $a^{t}$ defined in Section \ref{sec: defn}.
\begin{theorem}\label{thm 3}
 Let $G={\rm SO}(2n,\C)$. For each integer $0<k \leq n/2$, consider   $E_{n-k,k}=(C_{n-k,k},(C^{t}_{n-k,k})^{-1})$, where $C_{n-k,k}=(\underbrace{1,1,\cdots,1}_{n-k \ {\rm times}},\underbrace{-1,-1,\cdots,-1}_{k \ {\rm times}})$. Let $S_{[\underline{\lambda}]}\C^{2n}$ be the highest weight representation of ${\rm SO}(2n,\C)$ with highest weight $\underline{\lambda}=(\lambda_1,\lambda_2,\cdots,\lambda_n)$. Then:
 \begin{itemize}
     \item[A.]The character $\Theta_{\underline{\lambda}}(E_{n-k,k})=0$ if either $\#\eta_{0}(\underline{\lambda})>n-k$ or $\#\eta_{1}(\underline{\lambda})>n-k$. \\
     \item[B.]For $\#\eta_{0}(\underline{\lambda})=n-k$ or $\#\eta_{0}(\underline{\lambda})=k$, we have the following factorization:
\[
\Theta_{\underline{\lambda}}(E_{n-k,k})= 
\begin{cases}
	\pm 2^{e(k)} \cdot {\rm dim}(S_{[\lambda^{0}]}\C^{2n-2k} \otimes S_{[\lambda^{1}]}\C^{2k}) , & \text{if } \#\eta_{0}(\lambda)=n-k, \\
	 \pm 2^{e(k)} \cdot {\rm dim}(S_{[\lambda^{0}]}\C^{2k} \otimes S_{[\lambda^{1}]}\C^{2n-2k}) , & \text{if } \#\eta_{0}(\lambda)=k,	
\end{cases}
\]
where $e(k)=2\binom{n-2k}{2}-k+1$ and the highest weights $\lambda^{0}$ and $\lambda^{1}$ are defined by: \\ \\
$
\begin{cases}
	\begin{cases}
		\lambda^{0}+\overline{\rho}_{2n-2k}=\eta_{0}(\underline{\lambda})/2, \\
		\lambda^{1}+\overline{\rho}_{2k}=\eta_{1}(\underline{\lambda})/2.
	\end{cases} & \text{if} \ \#\eta_{0}(\lambda)=n-k, \\ \\
	\begin{cases}
		\lambda^{0}+\overline{\rho}_{2k}=\eta_{0}(\underline{\lambda})/2, \\
		\lambda^{1}+\overline{\rho}_{2n-2k}=\eta_{1}(\underline{\lambda})/2.
	\end{cases} & \text{if} \ \#\eta_{0}(\lambda)=k,
\end{cases}
$ \\ \\
where  $\overline{\rho}_{2n}=(n-1,n-2,\cdots,0)$ be the half the sum of the positive roots of ${\rm SO}(2n,\C)$. The highest weight representations $S_{[\lambda^{1}]}\C^{2k}$ and $S_{[\lambda^{1}]}\C^{2n-2k}$ are Spin representations of $\SO(2k,\C)$ and $\SO(2n-2k,\C)$ respectively.
 
 \item [C.] Let $\underline{\lambda}+\overline{\rho}_{2n}=(l_1,l_2,\cdots,l_n)$ and $n-k>k$. If $k<\#\eta_{0}(\lambda)<n-k$, we have the following alternating sum:
 \small
 \[
 \Theta_{\underline{\lambda}}(E_{n-k,k})=\pm \frac{\sum_{S \in I}\epsilon_{S}\cdot {\rm dim}(S_{ [ \lambda_{S} ]}\mathbb{C}^{2n-2k} \otimes S_{ [ \lambda_{S^{'}} ]}\mathbb{C}^{2k})}{2^{k(2n-2k-1)}},
 \]
where $S^{'}=[n] \backslash S$, $I=\{ \{i_1,i_2,\cdots,i_k\} \vert  \   l_{i_{j}  } \in \eta_{1}(\underline{\lambda})  \  \forall \ 1 \leq j \leq k\}$,  $\epsilon_{S}$ is the sign of the \enquote{shuffle} permutation corresponding to the sets of rows indexed by $S$ and $S^{'}$.
The highest weights $\lambda_{S}$ of  ${\rm SO}(2k,\C)$ and $\lambda_{S^{'}}$ of ${\rm SO}(2n-2k,\C)$ are given by:
\small
\begin{enumerate}
	\item[(i.)]$\lambda_{S}+\overline{\rho}_{2k}=(l_{s} \vert s \in S)$. 
	
	\item[(ii.)]$\lambda_{S^{'}}+\overline{\rho}_{2n-2k}=(l_{s} \vert s \in S^{'})$.
\end{enumerate}
\noindent There are a total of ${\# \eta_{1}(\underline{\lambda}) \choose k}$ terms in the above summation.
\end{itemize}

\begin{proof}
The proof of this theorem is analogous to the proof of the corresponding theorem for $\GL(n,\C)$ in which we replace the Weyl numerator (and the Weyl denominator) $\det(x^{\lambda_i+n-i}_{j})$ by $\det(x^{\lambda_i+n-i}_{j}-x^{-\lambda_i-n+i}_{j})+\det(x^{\lambda_i+n-i}_{j}+x^{-\lambda_i-n+i}_{j})$.  
\end{proof}
\end{theorem}    

In the next theorem we will discuss about the character theory of $\SO(2n+1,\C)$.
\section{The Odd Orthogonal Group ${\rm SO}(2n+1,\C)$}
For $a=(a_1,a_2,\cdots,a_{n}) \in (\C^{\ast})^{n}$, let us recall the notation $a^{t}$ defined in Section \ref{sec: defn}.
\begin{theorem}\label{thm 4}
 Let $G={\rm SO}(2n+1,\C)$. For each $0<k<n$, consider $F_{n-k,k}=(C_{n-k,k},1,(C^{t}_{n-k,k})^{-1})$, where $C_{n-k,k}=(\underbrace{1,1,\cdots,1}_{n-k \ {\rm times}},\underbrace{-1,-1,\cdots,-1}_{k \ {\rm times}})$. These are elements of order 2. Let $S_{[\underline{\lambda}]}\C^{2n}$ be the highest weight representation of ${\rm SO}(2n+1,\C)$ with highest weight $\underline{\lambda}=(\lambda_1,\lambda_2,\cdots,\lambda_n)$. 
 Let \(\underline{\lambda} + \rho_n = (l_1, l_2, \ldots, l_n)\). Then:      
     \begin{itemize}
     	\item[(A)] When $n=2k$, we have:
     	\[
     	\Theta_{\underline{\lambda}}(F_{k,k})=\pm {\rm dim}(S_{(\lambda^{0})}\C^{2k}),
     	\]
     	where the highest weight $\lambda^{0}$ of ${\rm GL}(2m,\C)$ is given by the receipe of Theorem 2.17 in \cite{AK}. \\
     	\item[(B.)] When $n-k \neq 
     	k$, the character $\Theta_{\underline{\lambda}}(F_{n-k,k})$ is given by: 
     	\[
     	\Theta_{\underline{\lambda}}(F_{n-k,k})= \frac{\sum_{L\subsetneq [n] \vert \#L=k}\epsilon_{L} \cdot {\rm dim}(S_{[\lambda_{L}]}\C^{2k}\bigotimes S_{[\lambda_{L^{'}}]}\C^{2n-2k+1})}{2^{2kn-2k^{2}+k-1}},
     	\]
     	where $[n]=\{1,2,\cdots,n\}$, $H=\{n-k+1,n-k+2,\cdots,n\}$, $L^{'}=[n] \textbackslash L$ and $\epsilon_{L}$ is the sign of the \enquote{shuffle} permutation corresponding to the sets of rows indexed by $L$ and $L^{'}$. 
The highest weights $\lambda_{L}$ of   $\SO(2k,\C)$ and $\lambda_{L^{'}}$ of $\SO(2n-2k+1,\C)$ are defined by: \\ \\
    \begin{enumerate}
     		\item[(a)]$\lambda_{L}+\overline{\rho}_{2k}=\{l_{u} \vert u \in L \}+\underbrace{(1/2,1/2,\cdots,1/2)}_{k \ {\rm times}}$.
     		\item[(b)]$\lambda_{L^{'}}+\overline{\rho}_{2n-2k+1}=\{l_{u} \vert u \in L^{'} \}+\underbrace{(1/2,1/2,\cdots,1/2)}_{n-k \ {\rm times}}$.
     	\end{enumerate} 
     \noindent For each $L \subset [n]$ with $\# L=k$, we have $S_{[\lambda_{L}]}\C^{2k}$ is a Spin representation of $\SO(2k,\C)$. 
     \end{itemize}
\begin{proof}
The proof for part (A) is done by Ayyer-Nishu in \cite{AK} and the proof for Part (B) is similar to the proof of Part (C) of Theorem \ref{thm 1}. So we will omit the proof of both Part (A) and Part (B) here.  
\end{proof}
\end{theorem}
\section{Comparisons between $A_{n}$,$B_{n}$,$C_{n}$ and $D_{n}$}

Let us focus on the classical groups \( C_4 = \mathrm{Sp}(8, \mathbb{C}) \) and \( B_4 = \mathrm{SO}(9, \mathbb{C}) \). Let \( \underline{\lambda} = (\lambda_1, \lambda_2, \lambda_3, \lambda_4) \) represent a highest weight. We analyze the highest weight representations \( S_{\langle \lambda \rangle} \mathbb{C}^8 \) associated with \( C_4 \) and \( S_{[\lambda]} \mathbb{C}^9 \) corresponding to \( B_4 \). Let us consider the element $D_{2,2}=(1,1,-1,-1,-1,-1,1,1)$ in $C_{4}$ and the element $F_{2,2}=(1,1,-1,-1,1,-1,-1,1,1)$ in $B_{4}$. 

Assuming that \( \lambda_1 \) is an even integer and \( \lambda_2, \lambda_3, \lambda_4 \) are odd integers, we invoke Theorem \ref{thm 2}. This theorem indicates that the character \( \Theta_{\underline{\lambda}}(D_{2,2}) \) vanishes, while the character \( \Theta_{\underline{\lambda}}(F_{2,2}) \) is non-zero. Specifically, the value of \( \Theta_{\underline{\lambda}}(F_{2,2}) \) can be expressed as:

\[
\Theta_{\underline{\lambda}}(F_{2,2}) = \pm \dim(S_{(\lambda^{1})} \mathbb{C}^4),
\]

where the highest weight \( \lambda^{1} \) for a representation of $\GL(4,\C)$ is defined by:

\[
\lambda^{1} + \rho_4 = \left( \frac{\lambda_1 + 4}{2}, \frac{\lambda_2 + 3}{2}, \frac{\lambda_4 + 1}{2}, \frac{-\lambda_3 - 1}{2} \right).
\]
This result follows either by a direct calculation or from \cite{AK} and underscores a crucial difference in the non-vanishing conditions for \( B_n \) compared to those for the classical groups \( A_n, C_n, \) and \( D_n \): the vanishing conditions for $A_{n},C_{n}$ and $D_{n}$ are the same but not for $B_{n}$. 
\begin{remark}
Unlike the characters of $\GL(n,\C)$, $\SO(2n,\C)$, $\Sp(2n,\C)$, there is no vanishing theorem about the character of the highest weight representation $S_{[\underline{\lambda}]}\C^{2n+1}$ with highest weight $\underline{\lambda}=(\lambda_1 \geq \lambda_2 \geq \cdots \geq \lambda_n \geq 0)$ at $F_{n-k,k}$.

\end{remark}
\section*{CHARACTER THEORY FOR $\mathrm{G}_{2}(\C)$}
The following proposition from the book of Fulton and Harris, Proposition 24.48 of \cite{FH} expressing the character of a highest weight representation of $\mathrm{G}_{2}(\mathbb{C})$ in terms of the character of the corresponding highest weight representation of $\mathrm{SL}_{3}(\mathbb{C}) \subset \mathrm{G}_{2}(\mathbb{C})$ is the main tool for us for computing the characters of $\mathrm{G}_{2}(\mathbb{C})$.

\begin{proposition}\label{prop 1}
 For the embedding of $\mathrm{SL}_{3}(\mathbb{C}) \subset \mathrm{G}_{2}(\mathbb{C})$, the character of the highest weight representation $\Pi_{k, l}$ of $\mathrm{G}_{2}(\mathbb{C})$ with highest weight $k \omega_{1}+l \omega_{2}$ is given by:
$$
\Theta_{k, l}=\frac{\mathcal{S}_{(k+2 l+1, k+l+1,0)}-\mathcal{S}_{(k+2 l+1, l, 0)}}{\mathcal{S}_{(1,1,0)}-\mathcal{S}_{(1,0,0)}},
$$
where $\mathcal{S}_{(a, b, c)}$ denotes the character of the highest weight representation of $\mathrm{GL}(3, \mathbb{C})$ with highest weight $(a, b, c)$ (restricted to $\mathrm{SL}_{3}(\mathbb{C})$ which is contained in $\mathrm{G}_{2}(\mathbb{C})$). Note that the highest weight of the representation of $\mathrm{SL}_{3}(\mathbb{C})$ with character $\mathcal{S}_{(k+2 l+1, k+l+1,0)}$ is the highest weight of the representation $\Pi_{k, l}$ of $\mathrm{G}_{2}(\mathbb{C})$ treated as a representation of the subgroup $\mathrm{SL}_{3}(\mathbb{C})$, and the character $\mathcal{S}_{(k+2 l+1, l, 0)}$ is that of the dual representation of $\mathrm{SL}_{3}(\mathbb{C})$; the denominator is the difference of the character of the standard 3 dimensional representation of $\mathrm{SL}_{3}(\mathbb{C})$ and its dual.
\end{proposition}
Proposition \eqref{prop 1} can be used to calculate the character of a highest weight representation of $\mathrm{G}_{2}(\mathbb{C})$ at any conjugacy class inside $\mathrm{SL}_{3}(\mathbb{C})$, hence at any conjugacy class of $\mathrm{G}_{2}(\mathbb{C})$. Since a calculation with the Weyl character formula for the conjugacy class $\mathcal{C}_{2}$ in $\mathrm{G}_{2}$, even with a more convenient variant such as Proposition 5.1, involves dealing with $0 / 0$, we have another result where calculation at the singular conjugacy class $\mathcal{C}_{2}$ is possible.

The following theorem when especialized to $x=1$ will calculate the character of irreducible representations of $\mathrm{G}_{2}(\mathbb{C})$ at $\mathcal{C}_{2}$. This theorem is at the same time, of theoretical importance in giving a factorization theorem just as in \cite{DP}, \cite{AK}.

\begin{theorem}\label{thm 5}
Let $\Pi_{k, l}$ be the highest weight representation of the group $\mathrm{G}_{2}(\C)$ with highest weight $k \omega_{1}+l \omega_{2}$. Then for an element $(a, b, c) \in \mathrm{SL}_{3}(\mathbb{C}) \subset \mathrm{G}_{2}(\mathbb{C})$ with $a \cdot b \cdot c=1$, we have the following character relationship at the particular element $\left(x,-x,-x^{-2}\right) \in$ $\mathrm{SL}_{3}(\mathbb{C}) \subset \mathrm{G}_{2}(\mathbb{C})$ \\
$\Theta_{k, l}\left(x,-x,-x^{-2}\right)=\left\{\begin{array}{lc}0, & \text { if } k, l \text { are odd, } \\ [1ex]
	-\Theta_{1}\left(x^{2}, x^{-2}\right) \times \Theta_{2}\left(x^{3}, x^{-3}\right), & \text { if } k \text { even }, l \text { odd, } \\ [1ex]
	\Theta_{3}\left(x^{2}, x^{-2}\right) \times \Theta_{4}\left(x^{3}, x^{-3}\right), & \text { if } k \text { odd }, l \text { even, } \\ [1ex] 
	-\Theta_{5}\left(x^{2}, x^{-2}\right) \times \Theta_{6}\left(x^{3}, x^{-3}\right), & \text { if } k, l \text { are even, }\end{array}\right.$ \\ \\
	here $\Theta_{1}, \Theta_{2}, \Theta_{3}, \Theta_{4}, \Theta_{5}$ and $\Theta_{6}$ are the characters of highest weight representations of $\mathrm{SL}_{2}(\mathbb{C})$ with highest weights $\frac{3 l+k}{4} L_{1}, \frac{k+l}{2} L_{1}, \frac{k-3}{4} L_{1}, \frac{k+2 l+1}{2} L_{1}, \frac{2 k+3 l+1}{4} L_{1}$, and $\frac{l-1}{2} L_{1}$, respectively.
\begin{proof}
We have
\[
\Theta_{k,l} = \frac{\mathcal{S}(k+2l+1, k+l+1, 0) - \mathcal{S}(k+2l+1, l, 0)}{\mathcal{S}(1, 1, 0) - \mathcal{S}(1, 0, 0)},
\]
where \( \mathcal{S}_{\lambda} \) denotes the Schur polynomial corresponding to the highest weight \( \lambda = (l_1, l_2, l_3) \), where \( l_1 \geq l_2 \geq l_3 \geq 0 \) are integers.

The Weyl denominator is easy enough to calculate:
\[
\mathcal{S}(1, 1, 0)(x, -x, -x^{-2}) - \mathcal{S}(1, 0, 0)(x, -x, -x^{-2}) = -\frac{x^2 - 1}{x^{2}}.
\]
Now we will calculate the Weyl numerator for different cases.

\textbf{Case I:} Suppose that \(k\) and \(l\) are both odd. In this case, the Schur polynomial \(\mathcal{S}(k+2l+1, k+l+1, 0)(x, -x, -x^{-2})\) can be expressed as:

\begin{align*}
\mathcal{S}(k+2l+1, k+l+1, 0)(x, -x, -x^{-2}) &= 
\frac{
	\begin{vmatrix}
		x^{k+2l+3} & (-x)^{k+2l+3} & (-x^{-2})^{k+2l+3} \\
		x^{k+l+2} & (-x)^{k+l+2} & (-x^{-2})^{k+l+2} \\
		1 & 1 & 1
	\end{vmatrix}
}{\Delta(x, -x, -x^{-2})} \\
 &= 0.
\end{align*}

Similarly, 
\begin{align*}
\mathcal{S}(k+2l+1, l, 0)(x, -x, -x^{-2}) &= 
\frac{
	\begin{vmatrix}
		x^{k+2l+3} & (-x)^{k+2l+3} & (-x^{-2})^{k+2l+3} \\
		x^{l+1} & (-x)^{l+1} & (-x^{-2})^{l+1} \\
		1 & 1 & 1
	\end{vmatrix}
}{\Delta(x, -x, -x^{-2})} \\
&=0.
\end{align*}
\textbf{Case II:} Suppose that \(k\) is even and \(l\) is odd. In this case, for \(X = (x, -x, -x^{-2})\), we have:
\begin{align*}
\mathcal{S}_{(k+2l+1,k+l+1,0)}(X) &= 
\frac{\begin{vmatrix}
		x^{k+2l+3} & (-x)^{k+2l+3} & (-x^{-2})^{k+2l+3} \\
		x^{k+l+2} & (-x)^{k+l+2} & (-x^{-2})^{k+l+2} \\
		1 & 1 & 1 \\
\end{vmatrix}}{\Delta(x, -x, -x^{-2})}
= 
\frac{\begin{vmatrix}
		x^{k+2l+3} & -x^{k+2l+3} & -(x^{-2})^{k+2l+3} \\
		x^{k+l+2} & -x^{k+l+2} & -x^{-2(k+l+2)} \\
		1 & 1 & 1 \\
\end{vmatrix}}{-2x \left( x + \frac{1}{x^2} \right) \left( x - \frac{1}{x^2} \right)} \\[2em]
&= 
\frac{\begin{vmatrix}
		0 & -x^{k+2l+3} & -(x^{-2})^{k+2l+3} \\
		0 & -x^{k+l+2} & -x^{-2(k+l+2)} \\
		2 & 1 & 1 \\
\end{vmatrix}}{-2x \left( x + \frac{1}{x^2} \right) \left( x - \frac{1}{x^2} \right)} 
= -\frac{\left( x^{-k-1} - x^{-k-3l-4} \right)}{\left( x^3 - \frac{1}{x^3} \right)}.
\end{align*}

\begin{align*}
	\mathcal{S}_{(k+2l+1, l, 0)}(X) &= 
	\frac{\begin{vmatrix}
			x^{k+2l+3} & (-x)^{k+2l+3} & (-x)^{k+2l+3} \\
			x^{l+1} & (-x)^{l+1} & (-x^{-2})^{l+1} \\
			1 & 1 & 1 \\
	\end{vmatrix}}{\Delta(x, -x, -x^{-2})} 
	= 
	\frac{\begin{vmatrix}
			x^{k+2l+3} & -x^{k+2l+3} & -x^{k+2l+3} \\
			x^{l+1} & -x^{l+1} & -x^{-2(l+1)} \\
			1 & 1 & 1 \\
	\end{vmatrix}}{-2x \left( x + \frac{1}{x^2} \right) \left( x - \frac{1}{x^2} \right)} \\[2em]
	&= 
	\frac{\begin{vmatrix}
			2x^{k+2l+3} & -x^{k+2l+3} & -x^{k+2l+3} \\
			0 & x^{l+1} & -x^{-2(l+1)} \\
			0 & 1 & 1 \\
	\end{vmatrix}}{-2x \left( x + \frac{1}{x^2} \right) \left( x - \frac{1}{x^2} \right)} 
	= -\frac{\left( x^{k+3l+4} - x^{k+1} \right)}{\left( x^3 - \frac{1}{x^3} \right)}.
\end{align*}
Therefore we have

\begin{align*}
	\mathcal{S}_{(k+2l+1,k+l+1,0)}(X) - \mathcal{S}_{(k+2l+1,l,0)}(X) 
	&= -\frac{\left( x^{-k-1} - x^{k+3l+4} \right) + \left( x^{k+1} - x^{-k-3l-4} \right)}{\left( x^3 - \frac{1}{x^3} \right)} \tag{5.2} \\[3em]
	&= -\frac{\left( \left( x^3 \right)^{\frac{l+1}{2}} - \left( x^{-3} \right)^{\frac{l+1}{2}} \right) \times \left( \left( x^2 \right)^{\frac{2k+3l+5}{4}} - \left( x^{-2} \right)^{\frac{2k+3l+5}{4}} \right)}{\left( x^3 - \frac{1}{x^3} \right)}. \tag{5.5}
\end{align*}

So we have
\begin{align*}
	\Theta_{k,l}(x, -x, -x^{-2})&= \frac{\left( (x^2)^{\frac{2k + 3l + 5}{4}} - (x^{-2})^{\frac{2k + 3l + 5}{4}} \right)}
	{(x^2 - x^{-2})} \times \frac{\left( (x^3)^{\frac{l + 1}{2}} - (x^{-3})^{\frac{l +1 }{2}} \right)}
	{(x^3 - x^{-3})} \\[10pt]
	&= \Theta_5(x^2, x^{-2}) \times \Theta_6(x^3, x^{-3}).
\end{align*}
where $\Theta_5$ and $\Theta_6$ are characters of the highest weight representations of $\mathrm{SL}_{2}(\mathbb{C})$ with highest weights 
$\frac{2k + 3l + 1}{4}$ and $\frac{l - 1}{2}$ respectively 
(note that these numbers may not themselves be integers when $k$ is even and $l$ odd, 
but their product is which is a reflection of the fact that we are not dealing with $\mathrm{SL}_2(\mathbb{C}) \times \mathrm{SL}_2(\mathbb{C})$ as a subgroup of $\mathrm{G}_2(\mathbb{C})$ 
but a quotient of it by $(-1, -1)$ which corresponds to $\mathrm{SL}_2(\alpha_2) \times \mathrm{SL}_2(\alpha_4)$ inside $\mathrm{G}_2(\C)$.

\textbf{Cases 3 and 4:} Both $k$ and $l$ are even, or $k$ is odd and $l$ is even. We omit a very similar calculation as before.
\end{proof}

\begin{remark}
It would be nice to understand what is behind the factorization of characters in Theorem \ref{thm 5}. We are not sure what is so special about the element $(x,-x,-x^{-2})$ used there for the factorization. For example at the element $(x,x,x^{-2})$, there is no factorization. 
\end{remark}

{\rm The proof of the following proposition is a direct consequence of Theorem 5.1 evaluated at the element $(x, -x, -x^{-2})$ for $x = 1$. This proposition is also available in the work of Reeder \cite{R}.}

\begin{proposition}\label{prop 2}
 Let $\Pi_{k,l}$ be the highest weight representation of $\mathrm{G}_2(\mathbb{C})$ with highest weight $k\omega_1 + l\omega_2$, where $\omega_1, \omega_2$ are the fundamental representations of $\mathrm{G}_2(\C)$, with character $\Theta_{k,l}(\mathcal{C}_2)$ at the unique conjugacy class $\mathcal{C}_2$ of order 2. Then we have
\[
\Theta_{k,l}(\mathcal{C}_2) = 
\begin{cases}
	0, & \text{if } k, l \text{ are odd} \\[15pt]
	\frac{(k + l + 2) \times (3l + k + 4)}{8}, & \text{if } k, l \text{ are even} \\[15pt]
	-\frac{(k + 1) \times (k + 2l + 3)}{8}, & \text{if } k \text{ odd, } l \text{ even} \\[15pt]
	-\frac{(l + 1) \times (3l + 2k + 5)}{8}, & \text{if } k \text{ even, } l \text{ odd}
\end{cases}
\]
\end{proposition}
\end{theorem}

\begin{remark}
From Proposition \ref{prop 2} it can be shown that when the character $\Theta_{k,l}(\mathcal{C}_2)$ is nonzero, up to a sign it can be written as half the dimension of a Spin representation of $\SO(4,\C)$. We have:

\[
\Theta_{k,l}(\mathcal{C}_2) = 
\begin{cases}
	\frac{1}{2} \ {\rm dim}(S_{[(\frac{k+2l+1}{2},\frac{l+1}{2})]}\C^{4}), & \text{if } k, l \text{ are even} \\[15pt]
	-\frac{1}{2} \ {\rm dim}(S_{[(\frac{k+l}{2},\frac{l+1}{2})]}\C^{4}), & \text{if } k \text{ odd, } l \text{ even} \\[15pt]
	-\frac{1}{2} \ {\rm dim}(S_{[(\frac{k+2l+1}{2},\frac{k+l+2}{2})]}\C^{4}), & \text{if } k \text{ even, } l \text{ odd}
\end{cases}
\]
\end{remark}
\noindent{\bf Acknowledgement:}
This work is part of the author's thesis work carried out at IIT Bombay under the guidance of Professor Dipendra Prasad whom the author thanks profusely for continued support on this work. 

\end{document}